\documentclass[11pt,dvips]{amsart}  
\usepackage{amssymb,amsmath,amsfonts}  
\usepackage[dvips]{graphicx,color}
\usepackage{cite}
\allowdisplaybreaks[4]
\definecolor{pgray}{gray}{0.8}

\newtheorem{theorem}{Theorem}[section]

\newtheorem{proposition}[theorem]{Proposition}
\newtheorem{lemma}[theorem]{Lemma}

\newtheorem{claim}[theorem]{Claim}

\numberwithin{equation}{section}

\title{\mbox{}}

\begin{document}
\begin{center}
{\bf \LARGE{
	Parabolic-scalings on large-time behavior of the incompressible Navier--Stokes flow\footnote{The completed version of this article is in Nonlinear Anal. Real World Appl. {\bf 85} (2025), 104350\\
	\hspace{6mm}\texttt{doi:10.1016/j.nonrwa.2025.104350}}
}}\\
\vspace{5mm}
Masakazu Yamamoto\footnote{e-mail : \texttt{masakazu@eng.niigata-u.ac.jp}\quad \texttt{mk-yamamoto@gunma-u.ac.jp}~ (current)}

\vspace{2mm}

Graduate School of Science and Technology, Niigata University\\
Ikarashi 2-no-cho 8050, Nishi-ku, Niigata 950-2181, Japan

\vspace{2mm}

Graduate School of Science and Technology, Gunma University (current)\\
Tenjin-cho 1-5-1, Kiryu 376-0052, Japan
\end{center}
\maketitle
\vspace{-15mm}
%
\begin{abstract}
Through asymptotic expansion, the large-time behavior of incompressible Navier--Stokes flow in $n$-dimensional whole space is depicted.
Especially, from their parabolic scalings, large-time behaviors of any terms on the expansion are clarified.
The parabolic scalings also guarantee the uniqueness of the expansion.
In the preceding work, the expansion with the $n$th order has already been derived.
They also predicted that the flow has some logarithmic evolutions in higher-order decay.
In this paper, an asymptotic expansion with $2n$th order is presented.
Furthermore, logarithmic evolutions are discovered.
\end{abstract}

\section{Introduction}
We study large-time behavior of the following initial-value problem of incompressible Navier--Stokes equations in whole space:
\begin{equation}\label{NS}
\left\{
\begin{array}{lr}
	\partial_t u + (u\cdot\nabla) u = \Delta u - \nabla p,
	&
	t>0,~ x \in \mathbb{R}^n,\\
	\nabla\cdot u = 0,
	&
	t>0,~ x \in \mathbb{R}^n,\\
	u(0,x) = a (x),
	&
	x \in \mathbb{R}^n,
\end{array}
\right.
\end{equation}
where $n \ge 2$, and $u = (u^1,u^2,\ldots,u^n) (t,x) \in \mathbb{R}^n$ and $p = p(t,x) \in \mathbb{R}$ are unknown velocity and pressure, respectively, and $a = (a^1,a^2,\ldots,a^n) (x) \in \mathbb{R}^n$ is given the initial velocity that satisfies $\nabla \cdot a = 0$.
The first and second equations come from the conservation laws of momentum and mass of fluid, respectively.
Under some conditions of smallness or smoothness of the initial data, $u$ exists globally in time and decays as $t \to +\infty$.
Particularly,
\begin{equation}\label{decayt}
	\| u(t) \|_{L^q (\mathbb{R}^n)} \le C t^{-\frac12} (1+t)^{-\gamma_q}
\end{equation}
is fulfilled for $1 \le q \le \infty$ and $\gamma_q = \frac{n}2 (1-\frac1q)$, and
\begin{equation}\label{decays}
	| u(t,x) | = O(|x|^{-n-1})
\end{equation}
as $|x| \to +\infty$ for any fixed $t$.
Those facts could be found in several works.
More precisely, we refer to \cite[Theorem 1]{Kt}, \cite[Theorem 4.1]{AmrchGrltSchnbk} and \cite[Theorem 4.2]{Gg-Mykw-Osd} for \eqref{decayt}, and \cite[Theorem 1.2]{Brndls-Vgnrn} and \cite[Theorem 1.1]{Mykw00} for \eqref{decays}.
We also refer to \cite{Brndls, Brndls-Okb, Cnnn-Krch, ChJn, FrwgKznShr08, FjtKt, Gg-Mykw, Kzn89JDE, KznOgwTnuch02} and the references therein for the historical background and various properties proven using \eqref{decayt} and \eqref{decays}.
For such a velocity, Carpio \cite{Crpo} and Fujigaki and Miyakawa \cite{Fjgk-Mykw} derived the asymptotic expansion of the Escobedo--Zuazua type \cite{EZ}.
Namely, they identified the terms $U_m \in C((0,\infty),L^1 (\mathbb{R}^n) \cap L^\infty (\mathbb{R}^n))$ which satisfy $\lambda^{n+m} U_m (\lambda^2 t, \lambda x) = U_m (t,x)$ for $\lambda > 0$, and
\begin{equation}\label{asymplow}
	\biggl\| u(t) - \sum_{m=1}^n U_m (t) \biggr\|_{L^q (\mathbb{R}^n)}
	\le
	C t^{-\gamma_q - \frac{n}2 - \frac12} \log (2+t)
\end{equation}
for $1 \le q \le \infty$ under the condition that $(1+|x|)^{n+1} a \in L^1 (\mathbb{R}^n)$.
Note that expansions with parabolic scalings are unique (for the details, see Subsection \ref{sbsect32}).
The goal is to extend this up to the $2n$th order and clearly present the logarithmic evolutions on the right-hand side of \eqref{asymplow}.
More precisely, we give unique functions $U_m$ and $K_m \in C((0,\infty),L^1 (\mathbb{R}^n) \cap L^\infty (\mathbb{R}^n))$ which satisfy $\lambda^{n+m} (U_m,K_m) (\lambda^2 t, \lambda x) = (U_m,K_m) (t,x)$ for $\lambda > 0$, and $u(t) \sim \sum_{m=1}^{2n} U_m (t) + \sum_{m=n+1}^{2n} K_m (t) \log t$ as $t \to +\infty$.
Note that this expansion contains some logarithmic evolutions.
It is well-known that the effect of viscosity does not provide such evolutions.
This characteristic evolution comes from the drift effect.
In other words, the higher-order expansion will capture the difference between Newtonian acceleration and Lagrangian acceleration.
In the expansion, some moments of solution are included as coefficients.
Hence, those moments should be integrable.
The slow decay \eqref{decays} prevents expansion to higher orders.
To achieve it, we use the knowledge of renormalization together with the idea via Kukavika and Reis \cite{Kkvc-Ris} (for the general theory of renormalization, see \cite{KtM,Ishg-Kwkm-Mchhs,Ishg-Kwkm,Ymd}).
They employed vorticity to derive spatial decay of the velocity.
The vorticity tensor $\omega^{ij} = \partial_i u^j - \partial_j u^i$ satisfies
\begin{equation}\label{v}
	\partial_t \omega^{ij} - \Delta \omega^{ij} +  \partial_i \mathcal{I}^j [u] - \partial_j \mathcal{I}^i [u] = 0,
\end{equation}
where
\[
	\mathcal{I}^j [u] = \omega^{\star j} \cdot u
\]
for $\omega^{\star j} = (\omega^{1j},\omega^{2j},\ldots,\omega^{nj})$.
Therefore a coupling of Duhamel principle and Biot--Savart law
\begin{equation}\label{BS}
	u^j = - \nabla (-\Delta)^{-1} \cdot \omega^{\star j}
\end{equation}
yields the mild solution of $u$ that
\begin{equation}\label{MSu}
\begin{split}
	u^j (t)
	&=
	- \nabla (-\Delta)^{-1} G(t) * \omega_0^{\star j}
	-
	\int_0^t
		\mathcal{R}^j \mathcal{R} G (t-s) * \mathcal{I} [u] (s)
	ds
	-
	\int_0^t
		G(t-s) * \mathcal{I}^j [u] (s)
	ds,
\end{split}
\end{equation}
where $\omega_0^{\star j} =  (\omega_0^{1j},\omega_0^{2j},\ldots,\omega_0^{nj})$ for $\omega_0^{ij} = \partial_i a^j - \partial_j a^i,~ \mathcal{R}^j = \partial_j (-\Delta)^{-1/2}$ and $\mathcal{R} = (\mathcal{R}^1,\mathcal{R}^2,\ldots,\mathcal{R}^n)$ are Riesz transforms, and $\mathcal{I} = (\mathcal{I}^1,\mathcal{I}^2,\ldots,\mathcal{I}^n)$ for the above $\mathcal{I}^j$.
Since $\omega$ is the vortex of incompressible fluid, it is natural that $\int_{\mathbb{R}^n} x^\alpha \omega_0 dx = 0$ for $|\alpha| \le 1$.
Indeed, Gauss' divergence theorem solves the case $\alpha = 0$, and the mass conservation law yields that $\int_{\mathbb{R}^n} x_h \omega_0^{ij} (x) dx = \int_{\mathbb{R}^n} (\delta_{hj} a^i - \delta_{hi} a^j) dx = 0$ for $1 \le h,i,j \le n$, where $\delta_{hi}$ is the Kronecker delta.
Due to the low frequency of $\nabla (-\Delta)^{-1}$ and $\mathcal{R}^j \mathcal{R}$, $u$ decays slowly as $|x| \to +\infty$ even if $\omega_0$ is localized.
Especially \eqref{decays} is optimal.
On the other hand, the decay rate of $\omega$ is controlled by $\omega_0$.
Kukavika and Reis exploited this fact to deal with the spatial structure of velocity.
Particularly, they proved that, upon some suitable conditions, $\omega$ fulfills that
\begin{equation}\label{decayog}
	\| \omega (t) \|_{L^q (\mathbb{R}^n)}
	\le
	C(1+t)^{-\gamma_q-1},\quad
	\| |x|^k \omega (t) \|_{L^q (\mathbb{R}^n)}
	\le
	C t^{-\gamma_q} (1+t)^{-1+\frac{k}2}
\end{equation}
for $1 \le q \le \infty,~ \gamma_q = \frac{n}2 (1-\frac1q)$ and some $k \in \mathbb{Z}_+$.
For such solutions, we publish the main results.
\begin{theorem}\label{thm1}
Let $\omega_0 \in L^1 (\mathbb{R}^n) \cap L^\infty (\mathbb{R}^n),~ |x|^{2n+1} \omega_0 \in L^1 (\mathbb{R}^n)$ and $\int_{\mathbb{R}^n} x^\alpha \omega_0$ $dx = 0$ for $|\alpha| \le 1$.
Assume that the solutions $u$ of \eqref{NS} for $a^j = - \nabla (-\Delta)^{-1} \cdot \omega_0^{\star j}$ and $\omega$ of \eqref{v} for $\omega (0) = \omega_0$ meet \eqref{decayt} and \eqref{decays} for $k = 2n+1$, respectively, where $\omega_0^{\star j}$ is the $j$th column.
Then there exist unique functions $U_m$ and $K_m \in C((0,\infty),L^1 (\mathbb{R}^n) \cap L^\infty (\mathbb{R}^n))$ such that
\begin{equation}\label{sc}
	\lambda^{n+m} (U_m,K_m) (\lambda^2 t, \lambda x) = (U_m,K_m) (t,x)
\end{equation}
for $\lambda > 0$, and
\begin{equation}\label{assertone}
	\biggl\| u(t) - \sum_{m=1}^{2n} U_m (t) - \sum_{m=n+1}^{2n} K_m (t) \log t\, \biggr\|_{L^q (\mathbb{R}^n)}
	= o(t^{-\gamma_q-n})
\end{equation}
as $t \to +\infty$ for $1 \le q \le \infty$.
In addition, if $|x|^{2n+2} \omega_0 \in L^1 (\mathbb{R}^n)$, then the left-hand side of \eqref{assertone} is estimated by $O(t^{-\gamma_q-n-\frac12} (\log t)^2)$
as $t \to +\infty$.
\end{theorem}
Among the terms included in the expansion, $U_{2n}$ decays fastest.
This term satisfies $\| U_{2n} (t) \|_{L^q(\mathbb{R}^n)} = t^{-\gamma_q-n} \| U_{2n} (1) \|_{L^q(\mathbb{R}^n)}$ because it has the parabolic scaling.
Therefore, the assertion \eqref{assertone} is sharp.
We emphasize that condition $|x|^{2n+1} \omega_0 \in L^1 (\mathbb{R}^n)$ is optimal because it is compatible with condition $|x|^{2n} a \in L^1 (\mathbb{R}^n)$.
The concrete forms of $U_m$ and $K_m$ will be discovered later.
Naturally, they satisfy $\nabla \cdot U_m = \nabla \cdot K_m = 0$ (see Appendix \ref{AppB}).
In the two-dimensional case, the expansion up to the fourth order has already been derived (cf. \cite{Ymmt}).
However, some terms in this expansion did not have parabolic scaling.
We solve the higher-dimensional cases and make some modifications in Claim \ref{clim31} in Section \ref{sect3}.
Now, in several-dimensional cases, we derive the terms that have the scalings to guarantee the uniqueness of the expansion (see the last sentences in Section \ref{sect3}).
Recently, the other equation on this motivation has been studied (cf. \cite{FkdIrn}).
They introduced a linearized problem and derived an asymptotic expansion and a logarithmic evolution from it.
This logarithmic evolution is similar to the term of $K_{n+1}$ in our assertion.

\vspace{3mm}

\paragraph{\bf Notations.}
For a vector $u$ and a tensor $\omega$, we denote their components $j$th and $ij$th by $u^j$ and $\omega^{ij}$, respectively.
We abbreviate the $j$th column of $\omega$ by $\omega^{\star j} = (\omega^{1j},\omega^{2j},\ldots,\omega^{nj})$.
For vector fields $f$ and $g$, the convolution of them is simply denoted by $f*g (x) = \int_{\mathbb{R}^n} f(x-y) \cdot g (y) dy = \int_{\mathbb{R}^n} f(y) \cdot g(x-y) dy$.
We often omit the spatial parameter from functions, for example, $u(t) = u(t,x)$.
In particular, $G(t) * \omega_0 = \int_{\mathbb{R}^n} G(t,x-y) \omega_0 (y) dy$ and $\int_0^t g(t-s) * f(s) ds = \int_0^t \int_{\mathbb{R}^n} g(t-s,x-y) f(s,y) dyds$.
We symbolize the derivations by $\partial_t = \partial/\partial t,~ \partial_j = \partial/\partial x_j$ for $1 \le j \le n,~ \nabla = (\partial_1,\partial_2,\ldots,\partial_n)$ and $\Delta = \lvert \nabla \rvert^2 = \partial_1^2 + \partial_2^2 + \cdots + \partial_n^2$.
The length of the multi-index $\alpha = (\alpha_1,\alpha_2,\ldots, \alpha_n) \in \mathbb{Z}_+^n$ is given by $\lvert \alpha \rvert = \alpha_1 + \alpha_2 + \cdots + \alpha_n$, where $\mathbb{Z}_+ = \mathbb{N} \cup \{ 0 \}$.
We abbreviate that $\alpha ! = \alpha_1 ! \alpha_2! \cdots \alpha_n !,~ x^\alpha = x_1^{\alpha_1} x_2^{\alpha_2} \cdots x_n^{\alpha_n}$ and $\nabla^\alpha = \partial_1^{\alpha_1} \partial_2^{\alpha_2} \cdots \partial_n^{\alpha_n}$.
We define the Fourier transform and its inverse by $\hat{\varphi} (\xi) = \mathcal{F} [\varphi] (\xi) = (2\pi)^{-n/2} \int_{\mathbb{R}^n} \varphi (x) e^{-ix\cdot\xi} dx$ and $\check{\varphi} (x) = \mathcal{F}^{-1} [\varphi] (x) = (2\pi)^{-n/2} \int_{\mathbb{R}^n} \varphi (\xi) e^{ix\cdot\xi} d\xi$, respectively, where $i = \sqrt{-1}$.
The Riesz transforms are defined by $\mathcal{R}^j \varphi = \partial_j (-\Delta)^{-1/2} \varphi = \mathcal{F}^{-1} [i\xi_j \lvert \xi \rvert^{-1} \hat{\varphi}]$ for $1 \le j \le n$ and $\mathcal{R} = (\mathcal{R}^1,\mathcal{R}^2,\ldots,\mathcal{R}^n)$.
Analogously, $\nabla (-\Delta)^{-1} \varphi = \mathcal{F}^{-1} [i\xi |\xi|^{-2} \hat{\varphi}]$.
The Lebesgue space and its norm are denoted by $L^q (\mathbb{R}^n)$ and $\| \cdot \|_{L^q (\mathbb{R}^n)}$, that is, $\| f \|_{L^q (\mathbb{R}^n)} = (\int_{\mathbb{R}^n} |f(x)|^q dx)^{1/q}$.
The heat kernel and its decay rate on $L^q (\mathbb{R}^n)$ are symbolized by $G(t,x) = (4\pi t)^{-n/2} e^{-|x|^2/(4t)}$ and $\gamma_q = \frac{n}2 (1-\frac1q)$.
We denote the floor function by the Gauss symbol $[\mu] = \max \{ m \in \mathbb{Z} \mid m \le \mu \}$.
We employ Landau's symbol.
Namely $f(t) = o(t^{-\mu})$ and $g(t) = O(t^{-\mu})$ mean $t^\mu f(t) \to 0$ and $t^\mu g(t) \to c$ for some $c \in \mathbb{R}$ such as $t \to +\infty$ or $t \to +0$.
Various positive constants are simply denoted by $C$.

\section{Preliminaries}
Before considering the velocity, we confirm the behavior of vorticity.
From \eqref{MSu}, the theory via Escobedo--Zuazua derive the asymptotic expansion of $u^j$ up to $n$th order by
\[
\begin{split}
	&U_m^j (t)
	=
	- \sum_{|\alpha| = m+1} \frac{\nabla^\alpha \nabla (-\Delta)^{-1} G(t)}{\alpha!} \cdot \int_{\mathbb{R}^n}
		(-y)^\alpha \omega_0^{\star j} (y)
	dy\\
	&- \sum_{2l+|\beta| = m} \frac{\partial_t^l \nabla^\beta \mathcal{R}^j \mathcal{R} G(t)}{l!\beta!} \cdot \int_0^\infty \int_{\mathbb{R}^n}
		(-s)^l (-y)^\beta \mathcal{I}[u] (s,y)
	dyds\\
	&- \sum_{2l+|\beta| = m} \frac{\partial_t^l \nabla^\beta G(t)}{l!\beta!} \int_0^\infty \int_{\mathbb{R}^n}
		(-s)^l (-y)^\beta \mathcal{I}^j [u] (s,y)
	dyds
\end{split}
\]
for $1 \le m \le n$.
Here $U_0 = 0$ since $\int_{\mathbb{R}^n} x^\alpha \omega_0 (x) dx= 0$ for $|\alpha| = 1$ and 
$
	\int_{\mathbb{R}^n} \mathcal{I}^i [u] dx = \int_{\mathbb{R}^n} (u\cdot\nabla u^i - u \cdot \partial_i u) dx = 0.
$
%
They satisfy $\lambda^{n+m} U_m (\lambda^2 t, \lambda x) = U_m (t,x)$ for $\lambda > 0$ and $\nabla \cdot U_m = 0$.
In fact, they are derived yet in \cite{Crpo,Fjgk-Mykw} essentially and satisfy \eqref{asymplow} (see Appendix \ref{SectA}).
Similarly, from mild solution of \eqref{v} for $\omega (0) = \omega_0$ (see \eqref{MSomg} in below), the theory yield the expansion also for $\omega^{ij}$ by
\begin{equation}\label{defOmg}
\begin{split}
	&\Omega_m^{ij} (t)
	=
	\sum_{|\alpha| = m} \frac{\nabla^\alpha G(t)}{\alpha!} \int_{\mathbb{R}^n} (-y)^\alpha \omega_0^{ij} (y) dy\\
	&+ \sum_{2l+|\beta|=m-1} \frac{\partial_t^l \nabla^\beta \partial_j G(t)}{l!\beta!} \int_0^\infty \int_{\mathbb{R}^n}
		(-s)^l (-y)^\beta \mathcal{I}^i [u] (s,y)
	dyds\\
	&- \sum_{2l+|\beta|=m-1}
	\frac{\partial_t^l \nabla^\beta \partial_i G(t)}{l!\beta!} \int_0^\infty \int_{\mathbb{R}^n}
		(-s)^l (-y)^\beta \mathcal{I}^j [u] (s,y)
	dyds
\end{split}
\end{equation}
for $2 \le m \le n+1$.
They also fulfill $\lambda^{n+m} \Omega_m (\lambda^2 t, \lambda x) = \Omega_m (t,x)$ for $\lambda > 0$.
It is not difficult to confirm the Biot--Savart law for them that $\Omega_{m}^{ij} = \partial_i U_{m-1}^j - \partial_j U_{m-1}^i$ and $U_m^j = - \nabla (-\Delta)^{-1} \cdot \Omega_{m+1}^{\star j}$, where $\Omega_{m+1}^{\star j}$ is the $j$th column (see also Appendix \ref{AppC}).
However, these relationships are never used.
The vorticity fulfills the following proposition.
\begin{proposition}\label{proppre}
Let $\omega_0 \in L^1 (\mathbb{R}^n) \cap L^\infty (\mathbb{R}^n),~ |x|^{n+2} \omega_0 \in L^1 (\mathbb{R}^n),~ \int_{\mathbb{R}^n} x^\alpha \omega_0$ $dx = 0$ for $|\alpha| \le 1$.
Assume that the solutions $u$ of \eqref{NS} for $a^j = - \nabla (-\Delta)^{-1} \cdot \omega_0^{\star j}$ and $\omega$ of \eqref{v} for $\omega (0) = \omega_0$ meet \eqref{decayt} and \eqref{decayog} for $k = n+2$, respectively, where $\omega_0^{\star j}$ is the $j$th column.
Then
\[
	\biggl\| \omega (t) - \sum_{m=2}^{n+1} \Omega_m (t) \biggr\|_{L^q (\mathbb{R}^n)}
	\le
	C t^{-\gamma_q-\frac{n}2-\frac12} (1+t)^{-\frac12} \log (2+t)
\]
for $1 \le q \le \infty$.
In addition, if $|x|^k \omega_0 \in L^1 (\mathbb{R}^n)$ and $\omega$ satisfies \eqref{decayog} for some $k \ge n+2$, then
\[
	\biggl\| |x|^k \biggl( \omega (t) - \sum_{m=2}^{n+1} \Omega_m (t) \biggr) \biggr\|_{L^q (\mathbb{R}^n)}
	\le
	C t^{-\gamma_q} (1+t)^{-\frac{n}2 - 1 + \frac{k}2} \log (2+t)
\]
for $1 \le q \le \infty$.
\end{proposition}
\begin{proof}
Duhamel principle provides the mild solution $\omega$ of \eqref{v} for $\omega (0) = \omega_0$ that
\begin{equation}\label{MSomg}
	\omega^{ij} (t) = G(t) * \omega_0^{ij} + \int_0^t \partial_i G(t-s) * \mathcal{I}^j [u] (s) ds - \int_0^t \partial_j G(t-s) * \mathcal{I}^i [u] (s) ds.
\end{equation}
Since the theory via Escobedo--Zuazua immediately provide it, we omit the proof of estimate without weight.
Another assertion is confirmed on the same way as in \cite{Ymmt}.
So reader may skip the following part.
The term of initial-data could be easily handled.
Indeed, since Taylor theorem yields that
\[
\begin{split}
	&|x|^k \biggl| G(t)*\omega_0^{ij} - \sum_{|\alpha|=2}^{n+1} \frac{\nabla^\alpha G(t)}{\alpha!} \int_{\mathbb{R}^n} (-y)^\alpha \omega_0^{ij} (y) dy \biggr|\\
	&\le
	C \sum_{|\alpha| = n+2} \int_{|y| \le |x|/2} \int_0^1 \lambda^{n+1} \left| |x-\lambda y|^k \nabla^\alpha G(t,x-\lambda y) (-y)^\alpha \omega_0^{ij} (y) \right| d\lambda dy\\
	&+ C \sum_{|\alpha|=2} \int_{|y| > |x|/2} \int_0^1 \lambda \left| \nabla^\alpha G(t,x-\lambda y) (-y)^\alpha |y|^k \omega_0^{ij} (y) \right| d\lambda dy\\
	&+ C \sum_{|\alpha|=2}^{n+1} \int_{|y| > |x|/2} \left| \nabla^\alpha G(t,x-y) (-y)^\alpha |y|^k \omega_0^{ij} (y) \right| d\lambda dy,
\end{split}
\]
Hausdorff-Young inequality and the scaling of $G$ complete the estimate.
Here we employed the assumption that $\int_{\mathbb{R}^n} (-x)^\alpha \omega_0^{ij} (y) dy = 0$ for $|\alpha| \le 1$.
We treat the external force terms on \eqref{MSomg}.
For them, we separate the domain $(0,\infty) \times \mathbb{R}^n$ to $Q_1 \cup Q_2 \cup Q_3 \cup Q_6$ or $Q_4 \cup Q_5 \cup Q_6$, where
\begin{gather*}
	Q_1 = (0,t/2] \times \{ y \in \mathbb{R}^n~ |~ |y| > |x|/2 \},\quad
	Q_2 = (0,t) \times \{ y \in \mathbb{R}^n~ |~ |y| \le |x|/2 \},\\
	Q_3 = (t/2,t) \times \{ y \in \mathbb{R}^n~ |~ |y| > |x|/2 \},\quad
	Q_4 = Q_2,\quad
	Q_5 = Q_1 \cup Q_3,\quad
	Q_6 = (t,\infty)\times \mathbb{R}^n.
\end{gather*}
Then
\[
\begin{split}
	&\int_0^t \partial_i G(t-s) * \mathcal{I}^j [u] (s) ds
	- \sum_{2l+|\beta|=1}^n \frac{\partial_t^l \nabla^\beta \partial_i G(t)}{l!\beta!} \int_0^\infty \int_{\mathbb{R}^n}
		(-s)^l (-y)^\beta \mathcal{I}^j [u] (s,y)
	dyds\\
	&=
	\rho_n^1 (t) + \rho_n^2 (t) + \cdots + \rho_n^6 (t),
\end{split}
\]
where
\[
\begin{split}
	&\rho_n^h (t)\\
	&=
	\left\{
	\begin{array}{lr}
	\displaystyle
		\iint_{Q_h}
			\biggl( \partial_i G(t-s,x-y) - \sum_{l=0}^{[n/2]} \frac{\partial_t^l \partial_i G(t,x-y)}{l!} (-s)^l \biggr)
	\mathcal{I}^j [u] (s,y)
		dyds,
		&
		h = 1,2,3,\\
	\displaystyle
		\sum_{l=0}^{[n/2]} \iint_{Q_h}
			\biggl( \frac{\partial_t^l \partial_i G(t,x-y)}{l!} - \sum_{|\beta|=0}^{n-2l} \frac{\partial_t^l \nabla^\beta \partial_i G(t,x)}{l!\beta!} (-y)^\beta \biggr)
		(-s)^l \mathcal{I}^j [u] (s,y)
		dyds,
		&
		h = 4,5,\\
	\displaystyle
		- \sum_{2l+|\beta| = 1}^n \frac{\partial_t^l \nabla^\beta \partial_i G(t)}{l!\beta!} \iint_{Q_h}
		(-s)^l (-y)^\beta \mathcal{I}^j [u] (s,y)
	dyds,
	&
		h=6.
	\end{array}
	\right.
\end{split}
\]
Taylor theorem also yields that
\[
\begin{split}
	|x|^k |\rho_n^1 (t)|
	&\le
	C \int_0^{t/2} \int_{|y| > |x|/2} \int_0^1 \lambda^{[\frac{n}2]} \left| \partial_t^{[\frac{n}2]+1} \partial_i G(t-\lambda s,x-y) \right|
	\left| (-s)^{[\frac{n}2]+1} |y|^k \mathcal{I}^j [u] (s,y) \right| d\lambda dyds,\\
	|x|^k |\rho_n^2 (t)|
	&\le
	C \int_0^{t/2} \int_{|y| \le |x|/2} \int_0^1 \lambda^{[\frac{n}2]} \left| |x-y|^k \partial_t^{[\frac{n}2]+1} \partial_i G(t-\lambda s,x-y)\right|
	\left| (-s)^{[\frac{n}2]+1} \mathcal{I}^j [u] (s,y) \right| d\lambda dyds
\end{split}
\]
and
\[
\begin{split}
	|x|^k |\rho_n^4 (t)|
	&\le
	C \sum_{2l+|\beta| = n+1} \int_0^t \int_{|y| \le |x|/2} \int_0^1
		\lambda^n \left| |x-\lambda y|^k \partial_t^l \nabla^\beta \partial_i G (t,x-\lambda y)\right|\\
    &\hspace{20mm}\left| (-s)^l (-y)^\beta \mathcal{I}^j [u] (s,y) \right|
	dyds.
\end{split}
\]
Hausdorff--Young inequality together with the scaling of $G$ and the properties of $\mathcal{I}$ give that
\[
\begin{split}
	&\| |x|^k \rho_n^1 (t) \|_{L^q (\mathbb{R}^n)} + \| |x|^k \rho_n^2 (t) \|_{L^q (\mathbb{R}^n)} + \| |x|^k \rho_n^4 (t) \|_{L^q (\mathbb{R}^n)}
    \le Ct^{-\gamma_q-\frac12} (1+t)^{-\frac{n}2-\frac12+\frac{k}2} \log (2+t).
\end{split}
\]
The remained terms $\rho_n^3$ and $\rho_n^6$ do not require Taylor theorem, in fact
\[
	\| |x|^k \rho_n^3 (t) \|_{L^q (\mathbb{R}^n)} + \| |x|^k \rho_n^6 (t) \|_{L^q (\mathbb{R}^n)} \le C t^{-\gamma_q-\frac12} (1+t)^{-\frac{n}2-\frac12+\frac{k}2}.
\]
The term $\rho_n^5$ is treated on the same way and some singularities as $t \to +0$ is mitigated by the scalings.
Thus, we conclude the proof.
\end{proof}
Those profiles satisfy $\lambda^{n+m} \Omega_m (\lambda^2 t, \lambda x) = \Omega_m (t,x)$ for $\lambda > 0$.
Summarizing the above fact, we see that $\mathcal{I} [u]$ in \eqref{MSu} has the products of $U_m$ and $\Omega_m$ as an asymptotic expansion.
Namely, from \eqref{asymplow} and Proposition \ref{proppre}, we see that
\begin{equation}\label{estI}
	\biggl\| |x|^{k} \biggl( \mathcal{I}[u] - \sum_{p=n+3}^{2n+2} \mathcal{I}_p \biggr) (t) \biggr\|_{L^q (\mathbb{R}^n)}
	\le
	C t^{-\gamma_q - n - 1 +\frac{k}2} (1+t)^{-\frac12} \log (2+t)
\end{equation}
for $1 \le q \le \infty$ and $0 \le k \le 2n+1$, where $\mathcal{I}_p = (\mathcal{I}_p^1,\mathcal{I}_p^2,\ldots,\mathcal{I}_p^n)$,
\begin{equation}\label{defIp}
	\mathcal{I}_p^j = \sum_{m=1}^{p-n-2} \Omega_{p-n-m}^{\star j} \cdot U_m
\end{equation}
and $\Omega_{p-n-m}^{\star j}$ is the $j$th column of \eqref{defOmg}.
By the way
\begin{equation}\label{scIp}
	\lambda^{n+p} \mathcal{I}_p (\lambda^2 t, \lambda x) = \mathcal{I}_p (t,x)
\end{equation}
for $\lambda > 0$ and
$
	\int_{\mathbb{R}^n} \mathcal{I}_p (1,x) dx = 0
$
are hold.
In fact, $\mathcal{I}_p^j = \sum_{m=1}^{p-n-m} \sum_{h=1}^n (\partial_h U_{p-n-m-1}^j - \partial_j U_{p-n-m-1}^h) U_m^h$ and $\sum_{h=1}^n \int_{\mathbb{R}^n} U_m^h \partial_h U_{p-n-m-1}^j dx = \sum_{m=1}^{p-n-m} \int_{\mathbb{R}^n} U_m^h \partial_j U_{p-n-m-1}^h dx = 0$ are valid.
Here, the second relation is confirmed by $\nabla \cdot U_m = 0$ and the elemental calculus.

In order to treat the velocity, we should study the integral kernels on \eqref{MSu}.
The basic properties of fractional integrations and the H\"ormander--Mikhlin estimates provide the following facts (cf. \cite{Sbt-Smz,Stin,Zmr}).
\begin{lemma}\label{lem-sc}
The Gaussian satisfies
\[
    \lambda^{n-1} \nabla (-\Delta)^{-1} G(\lambda^2 t, \lambda x) = \nabla (-\Delta)^{-1} G(t, x)
    \quad\text{and}\quad
    \lambda^n \mathcal{R}^j \mathcal{R} G(\lambda^2 t, \lambda x) = \mathcal{R}^j \mathcal{R} G(t,x)
\]
for $1 \le j \le n$ and $\lambda > 0$.
Moreover
\[
    (1+|x|)^{n-1+|\alpha|} \nabla^\alpha \nabla (-\Delta)^{-1} G(1,x)\quad\text{and}\quad
    (1+|x|)^{n+2l+|\alpha|} \partial_t^l \nabla^\alpha \mathcal{R}^j \mathcal{R} G(1,x)
\]
are in $L^\infty (\mathbb{R}^n)$ for $\alpha \in \mathbb{Z}_+^n$ and $l \in \mathbb{Z}_+$.
\end{lemma}
\begin{proof}
The parabolic scalings of the Gaussian are straightforward.
The spatial decays are coming from the H\"olmander--Mikhlin estimates since the singularities of the low frequencies of $\nabla^\alpha \nabla (-\Delta)^{-1} G$ and $\partial_t^l \nabla^\alpha \mathcal{R}^j \mathcal{R} G$ are estimated by $\xi^\alpha \xi |\xi|^{-2}$ and $|\xi|^{2l} \xi^\alpha \xi_j \xi |\xi|^{-2}$, respectively.
\end{proof}
\section{Proof of main results}\label{sect3}
Using the prepared tools, we show our main theorem.
Firstly, we introduce the rough expansion from \eqref{MSu} using the Taylor theorem together with the idea of renormalization.
Next, we adjust the shape of the expansion terms in Caims \ref{clim31} and \ref{clim32} to have parabolic scalings.
In particular, in Calim \ref{clim32}, we extract the logarithmic evolution from this expansion.
Claim \ref{propJ} confirms the integrability of the sensitive elements.
In these procedures, the scaling of $\mathcal{I}_p$ plays a crucial role.
In parallel, we estimate the remainder terms using the renormalization with the Carpio \cite{Crpo} and Fujigaki and Miyakawa \cite{Fjgk-Mykw} results.
Finally, we confirm the uniqueness of the expansion.
\subsection{Derivation of the expansion}
Now we expand the right-hand side of \eqref{MSu} and find $U_m$ and $K_m$ on \eqref{assertone}.
For some positive function $R(t) = o(t^{1/2})$ as $t \to +\infty$, the first term is expanded by Taylor theorem to
\[
\begin{split}
	&\nabla (-\Delta)^{-1} G (t) * \omega_0^{\star j}
	=
	\sum_{|\alpha| = 2}^{2n+1}
	\frac{\nabla^\alpha \nabla (-\Delta)^{-1} G (t)}{\alpha!}
	\cdot \int_{\mathbb{R}^n}
		(-y)^\alpha \omega_0^{\star j} (y)
	dy\\
	&+
	\sum_{|\alpha| = 2n+2} \int_{|y| \le R(t)} \int_0^1
		\lambda^{2n+1} \frac{\nabla^\alpha \nabla (-\Delta)^{-1} G (t,x-\lambda y)}{\alpha!} \cdot (-y)^\alpha \omega_0^{\star j} (y)
	d\lambda dy\\
	&+
	\sum_{|\alpha| = 2n+1} \int_{|y| > R(t)} \biggl( \int_0^1
		\lambda^{2n} \frac{\nabla^\alpha \nabla (-\Delta)^{-1} G (t,x-\lambda y)}{\alpha!} d\lambda
		- \frac{\nabla^\alpha \nabla (-\Delta)^{-1} G (t,x)}{\alpha!} \biggr)
		\cdot (-y)^\alpha \omega_0^{\star j} (y)
	dy.
\end{split}
\]
The first part provides some terms for $U_m$ since the scalings are guaranteed in Lemma \ref{lem-sc}.
The other parts decay faster than $t^{-\gamma_q-n}$ as $t\to +\infty$ in $L^q (\mathbb{R}^n)$ since $x^\alpha \omega_0 \in L^1 (\mathbb{R}^n)$ for $|\alpha| \le 2n+1$.
If $|x|^{2n+2} \omega_0 \in L^1 (\mathbb{R}^n)$ is assumed, then those other parts are bounded by $Ct^{-\gamma_q-n-\frac12}$ in $L^q (\mathbb{R}^n)$.
We expand the second term on \eqref{MSu}.
Firstly, we see the expansion as in \cite[Section 4]{Crpo} and \cite[Sections 4 and 5]{Fjgk-Mykw}.
More precisely, the Taylor expansion for the integral kernel on nonlinear term provides that
\begin{equation}\label{bsp}
\begin{split}
	&\int_0^t \mathcal{R}^j \mathcal{R}G (t-s) * \mathcal{I}[u] (s) ds
	=
	\sum_{2l+|\beta| = 1}^n \frac{\partial_t^l \nabla^\beta \mathcal{R}^j \mathcal{R}G(t)}{l!\beta!} \cdot \int_0^t \int_{\mathbb{R}^n}
		(-s)^l (-y)^\beta \mathcal{I}[u] (s,y)
	dyds\\
	&+
	\int_0^t \int_{\mathbb{R}^n} \biggl(
		\mathcal{R}^j \mathcal{R}G (t-s,x-y) - \sum_{2l+|\beta|=1}^n \frac{\partial_t^l \nabla^\beta \mathcal{R}^j \mathcal{R}G (t)}{l!\beta!} (-s)^l (-y)^\beta
	\biggr)
    \cdot \mathcal{I}[u] (s,y) dyds.
\end{split}
\end{equation}
The latter term of the right-hand side is separated to
\[
\begin{split}
	&\int_0^t \int_{\mathbb{R}^n} \biggl(
		\mathcal{R}^j \mathcal{R}G (t-s,x-y) - \sum_{2l+|\beta|=1}^n \frac{\partial_t^l \nabla^\beta \mathcal{R}^j \mathcal{R}G (t)}{l!\beta!} (-s)^l (-y)^\beta
	\biggr)
    \cdot \mathcal{I}[u] (s,y) dyds\\
	&=
	\sum_{2l+|\beta|=n+1} \frac{\partial_t^l \nabla^\beta \mathcal{R}^j \mathcal{R}G (t)}{l!\beta!} \cdot \int_0^t \int_{\mathbb{R}^n}
		(-s)^l (-y)^\beta \bigl( \mathcal{I}[u](s,y)
    - \mathcal{I}_{n+3} (1+s,y) \bigr)
	dyds\\
	&+
	\sum_{2l+|\beta|=n+1} \frac{\partial_t^l \nabla^\beta \mathcal{R}^j \mathcal{R}G (t)}{l!\beta!} \cdot \int_0^t \int_{\mathbb{R}^n}
		(-s)^l (-y)^\beta \mathcal{I}_{n+3} (1+s,y)
	dyds\\
	&+
	\int_0^t \int_{\mathbb{R}^n}
		\biggl( \mathcal{R}^j \mathcal{R}G (t-s,x-y) - \sum_{2l+|\beta|=1}^{n+1} \frac{\partial_t^l \nabla^\beta \mathcal{R}^j \mathcal{R}G (t,x)}{l!\beta!} (-s)^l (-y)^\beta \biggr)
		\cdot \mathcal{I}_{n+3} (s,y)
	dyds\\
	&+
	\int_0^t \int_{\mathbb{R}^n}
		\biggl( \mathcal{R}^j \mathcal{R}G (t-s,x-y) - \sum_{2l+|\beta|=1}^{n+1} \frac{\partial_t^l \nabla^\beta \mathcal{R}^j \mathcal{R}G (t,x)}{l!\beta!} (-s)^l (-y)^\beta \biggr)
		\cdot \left( \mathcal{I}[u] - \mathcal{I}_{n+3} \right) (s,y)
	dyds.
\end{split}
\]
We may expand the last term on the same way.
Hence we see that
\[
\begin{split}
	&\int_0^t \int_{\mathbb{R}^n}
		\biggl( \mathcal{R}^j \mathcal{R}G (t-s,x-y) - \sum_{2l+|\beta|=1}^{n+1} \frac{\partial_t^l \nabla^\beta \mathcal{R}^j \mathcal{R}G (t,x)}{l!\beta!} (-s)^l (-y)^\beta \biggr)
		\cdot \left( \mathcal{I}[u] - \mathcal{I}_{n+3} \right) (s,y)
	dyds\\
	&=
	\sum_{2l+|\beta|=n+2} \frac{\partial_t^l \nabla^\beta \mathcal{R}^j \mathcal{R} G(t)}{l!\beta!} \cdot \int_0^t \int_{\mathbb{R}^n}
		(-s)^l (-y)^\beta \bigl( (\mathcal{I}[u] - \mathcal{I}_{n+3}) (s,y)
    - \mathcal{I}_{n+4} (1+s,y) \bigr)
	dyds\\
	&+
	\sum_{2l+|\beta|=n+2} \frac{\partial_t^l \nabla^\beta \mathcal{R}^j \mathcal{R} G(t)}{l!\beta!} \cdot \int_0^t \int_{\mathbb{R}^n}
		(-s)^l (-y)^\beta \mathcal{I}_{n+4} (1+s,y)
	dyds\\
	&+
	\int_0^t \int_{\mathbb{R}^n}
		\biggl( \mathcal{R}^j \mathcal{R}G (t-s,x-y) - \sum_{2l+|\beta|=1}^{n+2} \frac{\partial_t^l \nabla^\beta \mathcal{R}^j \mathcal{R}G (t,x)}{l!\beta!} (-s)^l (-y)^\beta \biggr)
		\cdot \mathcal{I}_{n+4} (s,y)
	dyds\\
	&+
	\int_0^t \int_{\mathbb{R}^n}
		\biggl( \mathcal{R}^j \mathcal{R}G (t-s,x-y) - \sum_{2l+|\beta|=1}^{n+2} \frac{\partial_t^l \nabla^\beta \mathcal{R}^j \mathcal{R}G (t,x)}{l!\beta!} (-s)^l (-y)^\beta \biggr)\\
		&\hspace{40mm} \cdot \left( \mathcal{I}[u] - \mathcal{I}_{n+3} - \mathcal{I}_{n+4} \right) (s,y)
	dyds.
\end{split}
\]
By repeating this procedure, we get finally that
\begin{equation}\label{GI}
\begin{split}
	&\int_0^t \mathcal{R}^j \mathcal{R}G (t-s) * \mathcal{I}[u] (s) ds\\
	&=
	\sum_{2l+|\beta| = 1}^{2n} \frac{\partial_t^l \nabla^\beta \mathcal{R}^j \mathcal{R}G(t)}{l!\beta!} \cdot \int_0^t \int_{\mathbb{R}^n}
		(-s)^l (-y)^\beta\\
    &\hspace{20mm} \biggl( \Bigl( \mathcal{I}[u] - \sum_{p=2}^{2l+|\beta|+1} \mathcal{I}_p \Bigr) (s,y) - \mathcal{I}_{2l+|\beta|+2} (1+s,y) \biggr)
	dyds\\
	&+
	\sum_{2l+|\beta| = n+1}^{2n} \frac{\partial_t^l \nabla^\beta \mathcal{R}^j \mathcal{R}G(t)}{l!\beta!} \cdot \int_0^t \int_{\mathbb{R}^n}
		(-s)^l (-y)^\beta \mathcal{I}_{2l+|\beta|+2} (1+s,y)
	dyds\\
	&+
	\sum_{m=n+1}^{2n} \int_0^t \int_{\mathbb{R}^n}
		\biggl( \mathcal{R}^j \mathcal{R}G (t-s,x-y)\\
    &\hspace{20mm}
    - \sum_{2l+|\beta|=1}^m \frac{\partial_t^l \nabla^\beta \mathcal{R}^j \mathcal{R}G (t,x)}{l!\beta!} (-s)^l (-y)^\beta \biggr)
		\cdot \mathcal{I}_{m+2} (s,y)
	dyds
	+ r_{2n} (t),
\end{split}
\end{equation}
where
\[
\begin{split}
	r_{2n} (t)
    &=
	\int_0^t \int_{\mathbb{R}^n} \biggl(
		\mathcal{R}^j \mathcal{R}G (t-s,x-y) - \sum_{2l+|\beta| = 1}^{2n} \frac{\partial_t^l \nabla^\beta \mathcal{R}^j \mathcal{R}G (t,x)}{l!\beta!} (-s)^l (-y)^\beta
	\biggr)\\
	&\hspace{30mm}
	\cdot \Bigl( \mathcal{I}[u] - \sum_{p=n+3}^{2n+2} \mathcal{I}_p \Bigr) (s,y)
	dyds
\end{split}
\]
and we placed $\mathcal{I}_2 = \mathcal{I}_3 = \cdots = \mathcal{I}_{n+2} = 0$ for convenience.
The remainder term $r_{2n}$ satisfies
\[
	\| r_{2n} (t) \|_{L^q (\mathbb{R}^n)} \le C t^{-\gamma_q-n-\frac12} (\log t)^2
\]
as $t \to +\infty$ for $1 \le q \le \infty$.
Indeed, since Taylor theorem provides that
\[
\begin{split}
	&r_{2n} (t)\\
	&=
	\sum_{2l+|\beta| = 2n+1} \int_0^{t/2} \int_{\mathbb{R}^n} \int_0^1
		\lambda^n \frac{\partial_t^l \nabla^\beta \mathcal{R}^j \mathcal{R}G (t-\lambda s, x-\lambda y)}{l!\beta!}
		\cdot (-s)^l (-y)^\beta \Bigl( \mathcal{I}[u] - \sum_{p=n+3}^{2n+2} \mathcal{I}_p \Bigr) (s,y)
	d\lambda dyds\\
	&-
	\sum_{|\beta|=1} \int_{t/2}^t \int_{\mathbb{R}^n} \int_0^1
		(y \cdot \nabla) \mathcal{R}^j \mathcal{R}G (t-s,x-\lambda y)
		\cdot  \Bigl( \mathcal{I}[u] - \sum_{p=n+3}^{2n+2} \mathcal{I}_p \Bigr) (s,y)
	d\lambda dyds\\
	&-
	\sum_{2l+|\beta|=1}^{2n} \frac{\partial_t^l \nabla^\beta \mathcal{R}^j \mathcal{R}G (t)}{l!\beta!} \cdot \int_{t/2}^t \int_{\mathbb{R}^n}
		(-s)^l (-y)^\beta  \Bigl( \mathcal{I}[u] - \sum_{p=n+3}^{2n+2} \mathcal{I}_p \Bigr) (s,y)
	d\lambda dyds,
\end{split}
\]
Hausdorff--Young inequality and \eqref{estI} conclude the estimate.
We show in the following that the other parts in \eqref{GI} yield some terms in $U_m$ and $K_m$, and some error terms.
Now we recall Lemma \ref{lem-sc} on considerations of the coefficients on the first and second parts of \eqref{GI}.
\begin{claim}\label{clim31}
For $1 \le 2l+|\beta| \le 2n$ there exists a polynomial $\mathcal{P}$ of $(2n-2l-|\beta|)$th order such that
\[
\begin{split}
	&\int_0^t \int_{\mathbb{R}^n}
		(-s)^l (-y)^\beta \biggl( \Bigl( \mathcal{I}[u] - \sum_{p=2}^{2l+|\beta|+1} \mathcal{I}_p \Bigr) (s,y) - \mathcal{I}_{2l+|\beta|+2} (1+s,y) \biggr)
	dyds
	=
	\mathcal{P} (t^{-\frac12}) + o (t^{-n+l+\frac{|\beta|}2})
\end{split}
\]
as $t \to + \infty$, where $\mathcal{I}_p = 0$ for $2 \le p \le n+2$ and one for $n+3 \le p \le 2n+2$ is defined by \eqref{defIp}.
\end{claim}
\begin{proof}
We separate the interval on integration to $(0,t) = (0,\infty) \backslash [t,\infty)$ and correct $\mathcal{I}_{2l+|\beta|+2} (1+s)$ by $\mathcal{I}_{2l+|\beta|+2} (s)$ on $s \in [t,\infty)$, then
\[
\begin{split}
	&\int_0^t \int_{\mathbb{R}^n}
		(-s)^l (-y)^\beta \biggl( \Bigl( \mathcal{I}[u] - \sum_{p=2}^{2l+|\beta|+1} \mathcal{I}_p \Bigr) (s,y)- \mathcal{I}_{2l+|\beta|+2} (1+s,y) \biggr)
	dyds\\
	=
	&\int_0^\infty \int_{\mathbb{R}^n}
		(-s)^l (-y)^\beta \biggl( \Bigl( \mathcal{I}[u] - \sum_{p=2}^{2l+|\beta|+1} \mathcal{I}_p \Bigr) (s,y)- \mathcal{I}_{2l+|\beta|+2} (1+s,y) \biggr)
	dyds\\
	&-
	\int_t^\infty \int_{\mathbb{R}^n}
		(-s)^l (-y)^\beta \Bigl( \mathcal{I}[u] - \sum_{p=2}^{2l+|\beta|+2} \mathcal{I}_p \Bigr) (s,y)
	dyds\\
	&-
	\int_t^\infty \int_{\mathbb{R}^n}
		(-s)^l (-y)^\beta \left( \mathcal{I}_{2l+|\beta|+2} (1+s,y) - \mathcal{I}_{2l+|\beta|+2} (s,y) \right)
	dyds.
\end{split}
\]
On $[t,\infty)$, we expand the integrands to $\mathcal{I}[u] - \sum_{p=2}^{2l+|\beta|+2} \mathcal{I}_p \sim \sum_{p=2l+|\beta|+3}^{2n+2} \mathcal{I}_p$ and $\mathcal{I}_{2l+|\beta|+2} (1+s) - \mathcal{I}_{2l+|\beta|+2} (s) \sim \sum_{k = 1}^{[n-l-\frac{|\beta|}2]} \frac1{k!} \partial_s^k \mathcal{I}_{2l+|\beta|+2} (s)$.
Then we find that
\[
\begin{split}
	&\mathcal{P} (t^{-1/2})
	=
	\int_0^\infty \int_{\mathbb{R}^n}
		(-s)^l (-y)^\beta \biggl( \Bigl( \mathcal{I}[u] - \sum_{p=2}^{2l+|\beta|+1} \mathcal{I}_p \Bigr) (s,y) - \mathcal{I}_{2l+|\beta|+2} (1+s,y) \biggr)
	dyds\\
	&+
	\sum_{p=2l+|\beta|+3}^{2n+2}
	\int_t^\infty \int_{\mathbb{R}^n}
		(-s)^l (-y)^\beta \mathcal{I}_p (s,y)
	dyds
	+
	\sum_{k=1}^{[n-l-\frac{|\beta|}2]} \frac1{k!} \int_t^\infty \int_{\mathbb{R}^n}
		(-s)^l (-y)^\beta \partial_s^k \mathcal{I}_{2l+|\beta|+2} (s,y)
	dyds
\end{split}
\]
and
\[
\begin{split}
	&o(t^{-n+l+\frac{|\beta|}2})
	=
	- \int_t^\infty \int_{\mathbb{R}^n}
		(-s)^l (-y)^\beta \biggl( \mathcal{I}[u] - \sum_{p=n+3}^{2n+2} \mathcal{I}_p \biggr) (s,y)
	dyds\\
	&+
	\int_t^\infty \int_{\mathbb{R}^n}
		(-s)^l (-y)^\beta \biggl( \mathcal{I}_{2l+|\beta|+2} (1+s,y) - \sum_{k=0}^{[n-l-\frac{|\beta|}2]} \frac{\partial_s^k \mathcal{I}_{2l+|\beta|+2}(s,y)}{k!} \biggr)
	dyds.
\end{split}
\]
The decay of $\mathcal{I}[u]$ and the scaling of $\mathcal{I}_p$ say that $\mathcal{P}$ is sure a polynomial of $t^{-1/2}$ and the remained terms decay fast as $t \to +\infty$.
Indeed, the second and third parts of $\mathcal{P}$ are represented as
\[
	\int_t^\infty \int_{\mathbb{R}^n} (-s)^l (-y)^\beta \mathcal{I}_p (s,y) dyds = \frac{2 t^{-\frac{p}2 + l + \frac{|\beta|}2 + 1}}{p-2l-|\beta|-2} \int_{\mathbb{R}^n} (-1)^l (-y)^\beta \mathcal{I}_p (1,y) dy
\]
and
\[
\begin{split}
	&\int_t^\infty \int_{\mathbb{R}^n}
		(-s)^l (-y)^\beta \partial_s^k \mathcal{I}_{2l+|\beta|+2} (s,y)
	dyds
	=
	\frac{t^{-k}}k \int_{\mathbb{R}^n}
		(-1)^l (-y)^\beta \partial_s^k \mathcal{I}_{2l+|\beta|+2} (1,y)
	dy,
\end{split}
\]
since $p \ge 2l + |\beta| + 3$ and $k \ge 1$.
We should confirm that the first term of $\mathcal{P}$ is integrable even in $s \to +0$ and $s \to +\infty$.
Indeed, since $\| (-y)^\beta \mathcal{I}_p (s) \|_{L^1 (\mathbb{R}^n)} = s^{-\frac{p}2+\frac{|\beta|}2} \| (-y)^\beta \mathcal{I}_p (1) \|_{L^1 (\mathbb{R}^n)}$, we see that
\[
\begin{split}
	&\biggl| \int_{\mathbb{R}^n}
		(-s)^l (-y)^\beta \biggl( \Bigl( \mathcal{I}[u] - \sum_{p=2}^{2l+|\beta|+1} \mathcal{I}_p \Bigr) (s,y) - \mathcal{I}_{2l+|\beta|+2} (1+s,y) \biggr)
	dy \biggr|\\
	&\le
	C s^l \biggl( \| (-y)^\beta \mathcal{I}[u] (s) \|_{L^1 (\mathbb{R}^n)} + \sum_{p=2}^{2l+|\beta|+1} \| (-y)^\beta \mathcal{I}_p (s) \|_{L^1 (\mathbb{R}^n)}
    + \| (-y)^\beta \mathcal{I}_{2l+|\beta|+2} (1+s) \|_{L^1 (\mathbb{R}^n)} \biggr)\\
	&= O(s^{-\frac12})
\end{split}
\]
as $s \to +0$.
On the other hand, from \eqref{estI} and the scaling of $\mathcal{I}_p$, we have that
\[
\begin{split}
	&\biggl| \int_{\mathbb{R}^n}
		(-s)^l (-y)^\beta \biggl( \Bigl( \mathcal{I}[u] - \sum_{p=2}^{2l+|\beta|+1} \mathcal{I}_p \Bigr) (s,y) - \mathcal{I}_{2l+|\beta|+2} (1+s,y) \biggr)
	dy \biggr|\\
	&\le
	s^l \biggl\|
		(-y)^\beta \Bigl( \mathcal{I}[u] - \sum_{p=2}^{2l+|\beta|+2} \mathcal{I}_p \Bigr) (s)
	 \biggr\|_{L^1 (\mathbb{R}^n)}
	+
	s^l \biggl\| (-y)^\beta \left( \mathcal{I}_{2l+|\beta|+2} (1+s) - \mathcal{I}_{2l+|\beta|+2}(s) \right) \biggr\|_{L^1 (\mathbb{R}^n)}\\
	&= O(s^{-3/2} \log s)
\end{split}
\]
as $s \to +\infty$.
Therefore this is integrable in $s \in (0,\infty)$.
The fast decays of error terms are confirmed on the same way.
\end{proof}
%
Logarithmic evolutions, that is, the terms of $K_m$ come from the second part of \eqref{GI}.
\begin{claim}\label{clim32}
For $n+3 \le p \le 2n+2$ and $2l+|\beta|=p-2$ there exists a polynomial $\mathcal{Q}$ of $[n-l-\frac{|\beta|}2]$th order such that
\[
\begin{split}
	&\int_0^t \int_{\mathbb{R}^n}
		(-s)^l (-y)^\beta \mathcal{I}_p (1+s,y)
	dyds
	= \log t \int_{\mathbb{R}^n}
		(-1)^l (-y)^\beta \mathcal{I}_p (1,y)
	dy + \mathcal{Q} (t^{-1}) + o(t^{-n+l+\frac{|\beta|}2})
\end{split}
\]
as $t \to +\infty$.
\end{claim}
\begin{proof}
Since $2l+|\beta|=p-2$, the scaling of $\mathcal{I}_p$ immediately gives that
\[
\begin{split}
	&\int_0^t \int_{\mathbb{R}^n}
		(-s)^l (-y)^\beta \mathcal{I}_p (1+s,y)
	dyds
	=
	\int_0^t s^l (1+s)^{-l-1} ds \int_{\mathbb{R}^n} (-1)^l (-y)^\beta \mathcal{I}_p (1,y)
	dy.
\end{split}
\]
By separating $s^l = (1+s)^l + (s^l - (1+s)^l)$, we see
\[
    \int_0^t s^l (1+s)^{-l-1} ds
    =
    \log (1+t) + \sum_{k=0}^{l-1} \binom{l}{k} \int_0^t s^k (1+s)^{-l-1} ds.
\]
We repeat this separation under the integrations.
Finally, the elementary calculus and the local analyticity of the logarithmic function complete the proof.
\end{proof}
The last part on \eqref{GI} is well-defined and has their own parabolic scaling.
\begin{claim}\label{propJ}
For $n+1 \le m \le 2n$, the function
\[
\begin{split}
	J_m (t)
    &= \int_0^t \int_{\mathbb{R}^n}
		\biggl( \mathcal{R}^j \mathcal{R} G (t-s,x-y)
        - \sum_{2l+|\beta|=1}^m \frac{\partial_t^l \nabla^\beta \mathcal{R}^j \mathcal{R} G (t,x)}{l!\beta!} (-s)^l (-y)^\beta \biggr)
		\cdot \mathcal{I}_{m+2} (s,y)
	dyds
\end{split}
\]
is well-defined in $C((0,\infty),L^1 (\mathbb{R}^n) \cap L^\infty (\mathbb{R}^n))$ and has the parabolic scaling that $\lambda^{n+m} J_{m} (\lambda^2 t, \lambda x) = J_{m} (t,x)$ for $\lambda > 0$.
\end{claim}
\begin{proof}
By employing Taylor theorem, we have that
\[
\begin{split}
	&J_m (t)
	=
	\lim_{\varepsilon \to +0} \sum_{2l+|\beta|=m+1} \int_\varepsilon^{t/2} \int_{\mathbb{R}^n} \int_0^1
		\lambda^m \frac{\partial_t^l \nabla^\beta \mathcal{R}^j \mathcal{R} G (t-\lambda s, x-\lambda y)}{l!\beta!}
    \cdot (-s)^l (-y)^\beta \mathcal{I}_{m+2} (s,y)
	d\lambda dyds\\
	&-
	\int_{t/2}^t \int_{\mathbb{R}^n} \int_0^1 (y \cdot \nabla \mathcal{R}^j \mathcal{R} G) (t-s,x-\lambda y) \cdot \mathcal{I}_{m+2} (s,y) d\lambda dyds\\
	&- \sum_{2l+|\beta|=1}^m \frac{\partial_t^l \nabla^\beta \mathcal{R}^j \mathcal{R} G (t)}{l!\beta!} \cdot \int_{t/2}^t \int_{\mathbb{R}^n}
		(-s)^l (-y)^\beta \mathcal{I}_{m+2} (s,y)
	dyds.
\end{split}
\]
We care only the first term since integrability of the others are clear.
From the properties of $\mathcal{R}^j \mathcal{R} G$ and $\mathcal{I}_p$, we see that
\[
\begin{split}
	&\biggl\| \int_\varepsilon^{t/2} \int_{\mathbb{R}^n} \int_0^1
		\lambda^m \frac{\partial_t^l \nabla^\beta \mathcal{R}^j \mathcal{R} G (t-\lambda s, x-\lambda y)}{l!\beta!}
     \cdot (-s)^l (-y)^\beta \mathcal{I}_{m+2} (s,y)
	d\lambda dyds \biggr\|_{L^q (\mathbb{R}^n)}\\
	&\le
	C \int_{\varepsilon}^{t/2} \int_0^1 \lambda^m
		\| \partial_t^l \nabla^\beta \mathcal{R}^j \mathcal{R} G (t-\lambda s) \|_{L^q (\mathbb{R}^n)}
    \| (-s)^l (-y)^\beta \mathcal{I}_{m+2} (s) \|_{L^1 (\mathbb{R}^n)}
	d\lambda ds\\
	&\le
	C t^{-\gamma_q-l - \frac{|\beta|}2} \int_\varepsilon^{t/2} s^{-\frac{m}2-1+l+\frac{|\beta|}2} ds \| \mathcal{I}_{m+2} (1) \|_{L^1 (\mathbb{R}^n)}.
\end{split}
\]
Since $2l+|\beta|=m+1$, this converges as $\varepsilon \to +0$.
Therefore, from Lebesgue convergence theorem, the first term of $J_m$ also is integrable.
The parabolic scaling is confirmed by elementary calculations with Lemma \ref{lem-sc} and \eqref{scIp}.
\end{proof}
The last term on \eqref{MSu} could be treated by the same method.
By summarizing the above on \eqref{GI} and recalling Lemma \ref{lem-sc}, we conclude \eqref{assertone}.
\subsection{Uniqueness of the expansion}\label{sbsect32}
The uniqueness is straightforward.
Indeed, if \eqref{sc} and \eqref{assertone} are fulfilled for some $(\tilde{U}_m,\tilde{K}_m)$ instead of $(U_m,K_m)$, then
\[
	\| (U_1 - \tilde{U}_1)(t) \|_{L^q (\mathbb{R}^n)} \le \| (u - U_1)(t) \|_{L^q (\mathbb{R}^n)} + \| (u - \tilde{U}_1)(t) \|_{L^q (\mathbb{R}^n)} = o (t^{-\gamma_q-\frac12})
\]
as $t \to +\infty$ and
\[
	\| (U_1 - \tilde{U}_1)(t) \|_{L^q (\mathbb{R}^n)} = t^{-\gamma_q - \frac12} \| (U_1 - \tilde{U}_1)(1) \|_{L^q (\mathbb{R}^n)}.
\]
They are contradictory if $\tilde{U}_1 \neq U_1$.
The same procedure for $U_2 - \tilde{U}_2 = (u - \tilde{U}_1 - \tilde{U}_2) - (u-U_1-U_2)$ gives $\tilde{U}_2 = U_2$.
Inductively $\tilde{U}_m = U_m$ for $1 \le m \le n$.
Hence we see $\tilde{K}_{n+1} = K_{n+1}$ from 
\[
\begin{split}
	&\| (K_{n+1} - \tilde{K}_{n+1})(t) \log t \|_{L^q (\mathbb{R}^2)}\\
	&\le
	\biggl\| u(t) - \sum_{m=1}^n U_m (t) - K_{n+1} (t) \log t \biggr\|_{L^q (\mathbb{R}^n)}
	+
	\biggl\| u(t) - \sum_{m=1}^n \tilde{U}_m (t) - \tilde{K}_{n+1} (t) \log t \biggr\|_{L^q (\mathbb{R}^n)}\\
	&= o (t^{-\gamma_q-\frac{n}2-\frac12} \log t)
\end{split}
\]
as $t \to +\infty$ and
\[
	\| (K_{n+1} - \tilde{K}_{n+1})(t) \log t \|_{L^q (\mathbb{R}^2)} = t^{-\gamma_q - \frac{n}2 - \frac12} \log t \| (K_{n+1} - \tilde{K}_{n+1})(1) \|_{L^q (\mathbb{R}^2)}.
\]
Note that $\log t$ is not analytic in $t^{1/2}$ globally.
In this way, we finally confirm that $\tilde{U}_m = U_m$ and $\tilde{K}_m = K_m$ for $n+1 \le m \le 2n$, and then confirm our assertion.

\section*{Acknowledgments}
The author was supported by JSPS KAKENHI Grand Number 19K03560.
The author is very grateful to the reviewers for their careful reading and the valuable suggestions.

\appendix
\section{Approximation with low order is already well known}\label{SectA}
Usually, the mild solution of velocity is given by
\begin{equation}\label{HFK2}
	u^j (t) = G(t) * a^j - \sum_{k,h=1}^n \int_0^t (\delta_{jk} + \mathcal{R}^j \mathcal{R}^k) \partial_h G(t-s) * (u^h u^k) (s) ds,
\end{equation}
where $\delta_{jk}$ is the Kronecker delta.
As a natural conclusion, Carpio \cite[Theorem 0.6]{Crpo} and Fujigaki and Miyakawa \cite[Theorems 2.1 and 2.2]{Fjgk-Mykw} derived the asymptotic expansion up to $n$th order of it by
\[
\begin{split}
	&\tilde{U}_m^j (t)
	=
	\sum_{|\alpha| =m} \frac{\nabla^\alpha G(t)}{\alpha!} \int_{\mathbb{R}^n}
		(-y)^\alpha a^j (y)
	dy\\
	&- \sum_{2l+|\beta|=m-1} \sum_{k,h=1}^n \frac{(\delta_{jk} + \mathcal{R}^j \mathcal{R}^k) \partial_t^l \nabla^\beta \partial_h G(t)}{l!\beta!}
    \int_0^\infty \int_{\mathbb{R}^n}
		(-s)^l (-y)^\beta (u^h u^k) (s,y)
	dyds
\end{split}
\]
for $1 \le m \le n$.
In fact $\tilde{U}_m = U_m$, where $U_m$ was defined for our assertion.
Note that $\lambda^{n+m} \tilde{U}_m (\lambda^2 t, \lambda x) = \tilde{U}_m (t,x)$ for $\lambda > 0$.
Since the asymptotic expansion whose terms have parabolic scalings is unique (see the argument in Subsection \ref{sbsect32}), it is enough to confirm that \eqref{HFK2} and \eqref{MSu} are equivalent.
The Biot--Savart law $a^j = - \nabla (-\Delta)^{-1} \cdot \omega_0^{\star j}$ treats the terms of initial-data.
The same procedure as in \cite{Ymmt} estimates the nonlinear term.
Namely
\[
\begin{split}
	&\sum_{h=1}^n \int_0^t (\delta_{jk} + \mathcal{R}^j \mathcal{R}^k) \partial_h G(t-s) * (u^h u^k) (s) ds
	=
	\sum_{h=1}^n \int_0^t (\delta_{jk} + \mathcal{R}^j \mathcal{R}^k) G(t-s) * (u^h \partial_h u^k) (s) ds\\
	&=
	\sum_{h=1}^n \int_0^t (\delta_{jk} + \mathcal{R}^j \mathcal{R}^k) G(t-s) * (u^h \omega^{hk}) (s) ds
	+
	\sum_{h=1}^n \int_0^t (\delta_{jk} + \mathcal{R}^j \mathcal{R}^k) G(t-s) * (u^h \partial_k u^h) (s) ds
\end{split}
\]
and then
\[
\begin{split}
	&\sum_{k,h=1}^n \int_0^t (\delta_{jk} + \mathcal{R}^j \mathcal{R}^k) \partial_h G(t-s) * (u^h u^k) (s) ds\\
	&=
	\int_0^t
		\mathcal{R}^j \mathcal{R} G (t-s) * \mathcal{I} [u] (s)
	ds
	+
	\int_0^t
		G(t-s) * \mathcal{I}^j [u] (s)
	ds\\
	&+
	\frac12 \sum_{k=1}^n \int_0^t \mathcal{R}^j \mathcal{R}^k \partial_k G(t-s) * (|u|^2) (s) ds
	+
	\frac12 \int_0^t \partial_j G(t-s) * (|u|^2) (s) ds.
\end{split}
\]
The first and second terms on the right-hand side provide the nonlinear terms on \eqref{MSu}. The third and last terms are canceled since $\sum_{k=1}^n \mathcal{R}^j \mathcal{R}^k \partial_k = - \partial_j$.
This relation is valid in any suitable Sobolev spaces such as $W^{1,q} (\mathbb{R}^n)$, but it always holds at least for the Gaussian.

\section{Biot--Savart laws for the approximations}\label{AppC}
The relationships $U_m^j = - \nabla (-\Delta)^{-1} \cdot \Omega_{m+1}^{\star j}$ for $1 \le m \le n$ could be shown directly.
This is not difficult, but we select here the easiest method.
Hardy--Littlewood--Sobolev inequality (cf. \cite[Chapter V]{Stin} and \cite[Section 2]{Zmr}) yields that
\[
	\| \nabla (-\Delta)^{-1} \cdot (\omega^{\star j} - \Omega_2^{\star j}) \|_{L^{q_*} (\mathbb{R}^n)}
	\le
	C_q \| \omega^{\star j} - \Omega_2^{\star j} \|_{L^q (\mathbb{R}^n)}
\]
for $1 < q < n$ and $\frac1{q_*} = \frac1q - \frac1n$.
The right-hand side is estimated by $o (t^{-\gamma_{q_*}-\frac12})$ as $t \to +\infty$.
Hence, the Biot--Savart law \eqref{BS} says that $-\nabla (-\Delta)^{-1} \cdot \Omega_2^{\star j}$ is an approximation of $u^j$.
In particular, this function has the same scaling as $U_1^j$.
Hence, the uniqueness of the approximation provides that $U_1^j = -\nabla (-\Delta)^{-1} \cdot \Omega_2^{\star j}$.
The relationships for $2 \le m \le n$ are given inductively.

\section{Conservation-laws of the approximations}\label{AppB}
The terms $U_m$ and $K_m$ in Theorem \ref{thm1} satisfy $\nabla \cdot U_m = \nabla \cdot K_m = 0$.
Indeed, the same argument as in Subsection \ref{sbsect32} for $\nabla \cdot U_1 = \nabla \cdot (U_1 - u)$ together with $U_1^j - u^j = \nabla (-\Delta)^{-1} \cdot (\omega^{\star j} - \Omega_2^{\star j})$ says that $\nabla \cdot U_1 = 0$.
Here we used the bounds of the Riesz transforms in $L^q (\mathbb{R}^n)$ for some $1 < q < \infty$ (see \cite[Chapter III]{Stin} and also \cite[Section 2]{Zmr}) and the estimate for the vortex.
Inductively, we see the properties for any $1 \le m \le 2n$.
Using a similar procedure for $\int_{\mathbb{R}^n} x^\alpha \Omega_m dx = - \int_{\mathbb{R}^n} x^\alpha (\omega - \Omega_2 - \cdots - \Omega_m) dx$, we see that $\int_{\mathbb{R}^n} x^\alpha \Omega_m dx = 0$ for $|\alpha| \le 1$ and $2 \le m \le n+1$.


\begin{thebibliography}{99}
%
\bibitem{AmrchGrltSchnbk}
	Amrouche, C., Girault, V., Schonbek, M.E., Schonbek, T.P.,
	Pointwise decay of solutions and of higher derivatives to Navier--Stokes equations,
	SIAM J. Math. Anal. {\bf 31}  (2000), 740--753.
%
\bibitem{Brndls}
	Brandolese, L.,
	Space-time decay of Navier--Stokes flows invariant under rotations,
	Math. Ann. {\bf 329} (2004), 685--706.
%
\bibitem{Brndls-Vgnrn}
	Brandolese, L., Vigneron, F.,
	New asymptotic profiles of nonstationary solutions of the Navier--Stokes system,
	J. Math. Pures Appl. {\bf 88} (2007), 64--86.
%
%
\bibitem{Brndls-Okb}
	Brandolese, L., Okabe, T.,
	Annihilation of slowly-decaying terms of Navier--Stokes flows by external forcing,
	Nonlinearity {\bf 34} (2021), 1733--1757.
%
\bibitem{Cnnn-Krch}
	Cannone, M., Karch, G.,
	Smooth or singular solutions to the Navier--Stokes system?
	J. Differential Equations {\bf 197} (2004), 247--274.
%
\bibitem{Crpo}
	Carpio, A.,
	Large-time behavior in incompressible Navier--Stokes equation,
	SIAM J. Math. Anal. {\bf 27} (1996), 449--475.
%
%
\bibitem{ChJn}
	Choe, H.J, Jin, B.J.,
	Weighted estimate of the asymptotic profiles of the Navier--Stokes flow in $\mathbb{R}^n$,
	J. Math. Anal. Appl. {\bf 344} (2008), 353--366.
%
\bibitem{EZ}
	Escobedo, M., Zuazua, E.,
	Large time behavior for convection-diffusion equation in $\mathbb{R}^n$,
	J. Funct. Anal., {\bf 100} (1991), 119--161.
%
%
\bibitem{FrwgKznShr08}
	Farwig, R., Kozono, H., Sohr, H.,
	Criteria of local in time regularity of the Navier--Stokes equations beyond Serrin's condition,
	In: Parabolic and Navier--Stokes equations, Part 1, pp.175--184.
	Banach Center Publ., {\bf 81}, Part1,
	Polish Acad. Sci. Inst. Math., Warsaw. (2008)
%
\bibitem{Fjgk-Mykw}
	Fujigaki, Y., Miyakawa, T.,
	Asymptotic profiles of nonstationary incompressible Navier--Stokes flows in the whole space,
	SIAM J. Math. Anal. {\bf 33} (2001), 523--544.
%
\bibitem{FjtKt}
	Fujita, H., Kato, T.,
	On the Navier--Stokes initial value problem. I.,
	Arch. Rational Mech. Anal. {\bf 16} (1964), 269--315.
%
\bibitem{FkdIrn}
	Fukuda, I., Irino, Y.,
	Higher-order asymptotic profiles for solutions to the Cauchy problem for a dispersive-dissipative equation with a cubic nonlinearity,
	{\tt arXiv:2211.04667v2}.
%
\bibitem{Gg-Mykw}
	Giga, Y., Miyakawa, T.,
	Navier--Stokes flow in $\mathrm{R}^3$ with measures as initial vorticity and Morrey spaces,
	Comm. Partial Differential Equations {\bf 14} (1989), 577--618.
%
\bibitem{Gg-Mykw-Osd}
	Giga, Y., Miyakawa, T., Osada, H.,
	Two-dimensional Navier--Stokes flow with measures as initial vorticity,
	Arch. Rational Mech. Anal. {\bf 104} (1988), 223--250.
\bibitem{Ishg-Kwkm}
	Ishige, K., Kawakami, T.,
	Refined asymptotic expansions of solutions to fractional diffusion equations,
	J. Dynam. Differential Equations {\bf 36} (2024), 2679--2702.
	%
\bibitem{Ishg-Kwkm-Mchhs}
	Ishige, K., Kawakami, T., Michihisa, H.,
	Asymptotic expansions of solutions of fractional diffusion equations,
	SIAM J. Math. Anal. {\bf 49} (2017), 2167--2190.
%
%
\bibitem{KtM}
	Kato, M.,
	Sharp asymptotics for a parabolic system of chemotaxis in one space dimension,
	Differential Integral Equations {\bf 22} (2009), 35--51.
%
%
\bibitem{Kt}
	Kato, T.,
	Strong $L^p$-solutions of the Navier--Stokes equation in $\mathrm{R}^m$, with applications to weak solutions,
	Math. Z. {\bf 187} (1984), 471--480.
%
\bibitem{Kzn89JDE}
	Kozono, H.,
	Global $L^n$-solution and its decay property for the Navier--Stokes equations in half-space $\mathrm{R}_+^n$,
	J. Differential Equations {\bf 79} (1989), 79--88.
%
\bibitem{KznOgwTnuch02}
	Kozono, H., Ogawa, T., Taniuchi, Y.,
	The critical Sobolev inequalities in Besov spaces and regularity criterion to some semi-linear evolution equations,
	Math. Z. {\bf 242} (2002), 251-278.
%

%
\bibitem{Kkvc}
	Kukavica, I.,
	On the weighted decay for solutions of the Navier--Stokes system,
	Nonlinear Anal. {\bf 70} (2009), 2466-2470.
%
	%
\bibitem{Kkvc-Ris}
	Kukavica, I., Reis, E.,
	Asymptotic expansion for solutions of the Navier--Stokes equations with potential forces,
	J. Differential Equations {\bf 250} (2011), 607--622.
 %
%
\bibitem{Kkvc-Trrs}
	Kukavica, I., Torres, J.J.,
	Weighted $L^p$ decay for solutions of the Navier--Stokes equations,
	Comm. Partial Differential Equations {\bf 32} (2007), 819--831.
%
%
\bibitem{Mykw00}
    Miyakawa, T.,
    On space-time decay properties of nonstationary incompressible Navier--Stokes flows in $\mathrm{R}^n$,
    Funkcial. Ekvac. {\bf 43} (2000), 541--557.
%

%
\bibitem{Sbt-Smz}
	Shibata, Y., Shimizu S.,
	A decay property of the Fourier transform and its application to the Stokes problem,
	J. Math. Fluid Mech. {\bf 3} (2001), 213--230.
%
\bibitem{Stin}
	Stein, E.M.,
	Singular Integrals and Differentiability Properties of Functions,
	Princeton University Press,
	Princeton, New Jersey, 1970.
%
\bibitem{Ymd}
	Yamada, T.,
	Higher-order asymptotic expansions for a parabolic system modeling chemotaxis in the whole space,
	Hiroshima Math. J. {\bf 39} (2009), 363--420.
%
\bibitem{Ymmt}
	Yamamoto, M.,
	Time evolution of the Navier--Stokes flow in far-field,
	J. Math. Fluid Mech.  {\bf 26} (2024), 67.
%
\bibitem{Zmr}
	Ziemer, W.P.,
	Weakly Differentiable Functions,
	Graduate Texts in Math., vol. 120, Springer Verlag, New York, 1989.
\end{thebibliography}
\end{document}